\documentclass{jancl}

\journal{13}{1}{2003}{-}{-}

\usepackage[all]{xy}

\newjancltheorem{oq}{Open question}{\itshape}

\renewcommand{\a}{\alpha}
\renewcommand{\b}{\beta}
\newcommand{\g}{\gamma}

\newcommand{\X}{\chi}

\newcommand{\f}{\rightarrow}
\newcommand{\q}{\forall}
\newcommand{\faux}{\perp}
\newcommand{\non}{\neg}
\newcommand{\G}{\Gamma}
\newcommand{\D}{\Delta}

\renewcommand{\v}{\vdash}

\newcommand{\pp}{\leq}
\newcommand{\pg}{\geq}
\newcommand{\ou}{\vee}

\newcommand{\et}{\vedge}
\renewcommand{\int}{\bigcap}

\renewcommand{\et}{\wedge}

\newcommand{\inclus}{\subseteq}

\newcommand{\F}{\displaystyle\frac}

\newcommand{\realise}{\ifmmode{ \v\mkern-7.6mu
\v\mkern2mu}\else{ \v\kern-.28em
\v\hskip1pt\ }\fi\relax}

\title[Propositional Mixed Logic]{Propositional Mixed Logic:
Its Syntax and\\ Semantics}

\author{Karim Nour\fup{*} \andauthor\/ Abir Nour\fup{**}}

\address{%
Université de Chambéry \\
LAMA -- \'Equipe de Logique \\
73376 Le Bourget du Lac \\[3pt]
nour@univ-savoie.fr\\[6pt]
Université libanaise \\
Faculté des sciences \\
Section III - BP 826 \\
Tripoli - Liban \\[3pt]
abir\_n@inco.com.lb
}

\abstract{In this paper, we present a propositional logic (called mixed logic)
containing disjoint copies of minimal, intuitionistic and classical
logics. We prove a completeness theorem for this logic with respect to
a Kripke semantics. We establish some relations between mixed logic
and minimal, intuitionistic and classical logics. We present at the
end a sequent calculus version for this logic.}

\keywords{propositional logic, minimal logic, intuitionistic
logic, classical logic, mixed logic, completeness theorem, Kripke
semantics.}

\begin{document}

\maketitle

\section{Introduction}

Propositional intuitionistic and classical logics (abbreviated: PLI and
PLC) are built by adding absurdity rules to propositional minimal logic
(abbreviated PLM). The best known formalization consists to adding the
intuitionistic absurdity rule (from the absurdity we can deduce all
formulas) to PLM to obtain PLI, and to adding the classical absurdity rule
(a non false formula is true) to PLM (or PLI) to obtain PLC. With this
kind of formalism there are some problems.

\begin{itemize}
\item A classical formula does not contain any information on the
smallest logical system in which it is derivable. To have this
information, we must use the non effective decision algorithms of
PLM and PLI. But with these algorithms we cannot know how many
times we used the absurdity rules and on which formulas.

\item A formula has several derivations and the formula does not
contain informations to find its ``better'' derivation. For
example, if one takes $A = (X \f Y) \ou (Y \f X)$, we can prove
this formula using the classical absurdity rule on $A$ (i.e.\  we
prove $\non \non A$).  And we can also prove it using the
classical absurdity rule on the variable $Y$. Indeed, if $Y$ is
true, then we have (in PLM) $X \f Y$, and if $Y$ is false, then we
have (in PLI) $Y \f X$. The second derivation is nearer to the
human reasoning. For this reason we want to call it ``a good
derivation'' of the formula $A$.

\item Each of these three logics has a semantics and a
completeness theorem. For PLC it is the truth tables, for PLI it
is the intuitionistic Kripke models and for PLM it is the minimal
Kripke models. If we look closely at the proofs of the
completeness theorems, a great resemblance is seen. Why not study
all these logics at the same time? i.e.\  introduce a single
semantics for these logics and only prove one completeness theorem
in order to deduce the completeness of each system.
\end{itemize}

We propose in this paper a partial solution to these problems. We
present a propositional logic (called mixed logic and abbreviated PML)
containing three kinds of variables: minimal variables indexed by
$m$, intuitionistic variables indexed by $i$ and classical variables
indexed by $c$. We restrict the absurdity rules to the formulas
containing the corresponding variables. The main novelty of our system
is that minimal, intuitionistic and classical logics appear as
fragments. For instance a proof of an intuitionistic formula may use
classical lemmas without any restriction. This approach is radically
different from the one that consists in changing the rule of the game
when we want to change logic. Here there is only one logic which,
depending on its use, may appear classical, intuitionistic or
minimal. We introduce for the system PML a Kripke semantics which is the
superposition of minimal, intuitionistic and classical semantics. We
show a completeness theorem which implies the completeness theorems of
systems PLM, PLI and PLC. We deduce from this theorem a very
significant result which is the following: ``for a formula $A$ to be
derivable in a logic, it is necessary that the formula contains at
least a variable which corresponds to this system''. We were
interested by labelling problems (we label variables by $m$, $i$ or
$c$) for classical formulas. We present decision algorithms for these
problems and we formally define the concept of ``good
derivation'' for a classical formula. We also present a sequent
calculus version of this system. This presentation is coherent with
what we already know on sequent calculus: classical logic comes from
the possibility to put several formulas on the
right.

This paper is an introduction to this domain and much questions remain
open. For example, the standard proofs of cut-elimination are not
adapted to our system. This comes primarily from impossibility of
coding disjunction.

The idea to present only one system for different logics is not
completely new.  Indeed, J.Y.\ Girard presented in \cite{GIR:1993}
a single sequent calculus (denoted LU) common to classical,
intuitionistic and linear logics. The idea of Girard is to use a
single variable set but different connectives which correspond to
each fragment. Each formula is given with a polarity: positive,
neutral and negative. For each connective the rules depend on the
polarity of the formulas. On the other hand the system LU has a
cut-elimination theorem and then the sub-formula property.

Finally, let us mention that J.-L.\ Krivine and K.\ Nour
introduced a second order mixed logic in order to type storage and
control operators in $\lambda$-calculus (see \cite{NOU:2000}). The
theoretical properties of this system are not difficult to prove
because the only connectives are $\f$ and $\q$. The presence of
$\ou$ in system PML complicates our study.

\section{The system PML}

We present in this section the natural deduction version of
propositional mixed logic.

\begin{definitions}
\begin{janclenum}{(\arabic}{)}
\item We suppose that we have three disjoint countable sets of
propositional variables: ${\cal V}_m = \{X_m , Y_m , Z_m , ... \}$
the set of \emph{minimal variables}, ${\cal V}_i = \{X_i , Y_i ,
Z_i , ... \}$ the set of \emph{intuitionistic variables}, ${\cal
V}_c = \{X_c , Y_c , Z_c , ... \}$ the set of \emph{classical
variables} and a special constant denoted $\perp$.

\item The \emph{formulas} are defined by induction. Each element
of ${\cal P} = {\cal V}_m \cup {\cal V}_i \cup {\cal V}_c \cup \{
\perp \}$ is a formula. And if $A$,$B$ are formulas, then $A \et
B$, $A \ou B$ and $A \f B$ are formulas.  We denote  $\neg A = A
\f \perp$.

\item If $A$ is a formula, we denote by $var(A)$ the set of
variables of $A$.  A \emph{classical formula} (resp.\ an
\emph{intuitionistic
  formula}) is a formula $A$ such that $var(A) \subseteq {\cal V}_c$
(resp.\ $var(A) \subseteq {\cal V}_i \cup {\cal V}_c$). We allow
the use of classical variables to build intuitionistic formulas
because the intuitionistic absurdity rule is derivable in
classical logic.

\item A \emph{simple sequent} is an expression of the form $\G \v
A$ where $\G \cup \{ A \}$ is a finite set of formulas. A
derivation ${\cal D}$ may be constructed according to one of the
rules below.
\end{janclenum}

\begin{minipage}[t]{120pt}
$(Ax) \; \F{}{A \v A}$\\
\end{minipage}
\begin{minipage}[t]{130pt}
$(W) \; \F{\G \v A} {\G , B \v A}$\\
\end{minipage}

\begin{minipage}[t]{120pt}
$(\et_I) \; \F{\G_1 \v A_1 \;\;\; \G_2 \v A_2 }
{\G_1 , \G_2 \v A_1 \et A_2}$\\
\end{minipage}
\begin{minipage}[t]{130pt}
$(\et_E) \; \F{\G \v A_1 \et A_2 } { \G \v A_i}$ \\
\end{minipage}

\begin{minipage}[t]{120pt}
$(\ou_I) \;\F{\G \v A_i} {\G \v A_1 \ou A_2}$ \\
\end{minipage}
\begin{minipage}[t]{195pt}
$(\ou_E) \; \F{\G_1 \v A_1 \ou A_2  \;\;\; \G_2 , A_1 \v B  \;\;\; \G_3 , A_2 \v B}  { \G_1 , \G_2 , \G_3 \v B}$\\
\end{minipage}

\begin{minipage}[t]{120pt}
$(\f_I) \; \F{\G , A_1 \v A_2} { \G \v A_1 \f A_2}$\\
\end{minipage}
\begin{minipage}[t]{160pt}
$(\f_E) \; \F{\G_1 \v A_1 \f A_2 \;\;\;  \G_2 \v A_1} {\G_1 , \G_2 \v A_2 }$\\
\end{minipage}

\begin{minipage}[t]{185pt}
$(\perp_i) \;\F{\G \v \perp  \;\;\;   A {\rm\;  is \;  an \; intuitionistic \; formula }}{ \G\v A}$\\
\end{minipage}

\begin{minipage}[t]{180pt}
 $(\perp_c) \; \F{\G \v \neg \neg A \;\;\;   A {\rm\;  is \;  a \; classical \; formula }  }{ \G\v A}$
\end{minipage}

The rules given above determine the natural deduction system,
abbreviated \emph{PML}. If ${\cal D}$ is a derivation ending with
a simple sequent $\G \v A$, then we write $\G \v_{pml} A$.
\end{definitions}

\begin{example}\label{dem}
\begin{janclenum}{\alph}{)}
\item $\v_{pml}  X_c \ou \neg X_c$.

\begin{center}
{\footnotesize $\F{\F{\F{\F{\F{\F{\F{\F{}{X_c \v X_c}}{X_c \v X_c
\ou \neg X_c} \quad \F{}{\neg (X_c \ou \neg X_c) \v \neg (X_c \ou
\neg X_c)}} {X_c , \neg (X_c \ou \neg X_c) \v  \perp}}{\neg (X_c
\ou \neg X_c) \v
    \neg X_c}}
{\neg (X_c \ou \neg X_c) \v X_c \ou \neg X_c}
\quad
\F{}{\neg (X_c \ou \neg X_c) \v \neg (X_c \ou \neg X_c)}}
{\neg (X_c \ou \neg X_c) \v \perp}}
{\v \neg \neg (X_c \ou \neg X_c)}}
{\v  X_c \ou \neg X_c}$}
\end{center}

\item $\v_{pml} (X_m \f X_c) \ou (X_c \f X_i)$.

\begin{center}
{\footnotesize $\F{\F{\vdots}{\v X_c \ou \non X_c} \quad
\F{\F{\F{\F{}{X_c \v X_c}}{X_c , X_m \v X_c}}{X_c \v X_m \f
X_c}}{X_c
 \v (X_m \f X_c) \ou (X_c \f X_i)}
 \quad
\F{\F{\F{\F{\F{}{X_c \v X_c} \quad \F{}{\neg X_c \v \neg X_c}}{X_c , \neg X_c \v \perp}}{X_c , \neg X_c \v X_i}}{
\neg X_c \v  X_c \f X_i}}{\neg X_c \v (X_m \f X_c) \ou (X_c \f X_i)}}
{\v (X_m \f X_c) \ou (X_c \f X_i)}$}
\end{center}

\item  $\v_{pml} ( X_c \f X_m \ou X_i ) \f ( X_m \ou ( X_c \f X_i
))$ (left to the readers).
\end{janclenum}
\end{example}

\begin{remark} {Note that the indices of variables used in the derivable formulas give
some ideas on their derivations. For the formula $(X_m \f X_c) \ou
(X_c \f X_i)$, the classical absurdity rule is used on the variable
$X_c$ and the intuitionistic absurdity rule is used on the variable
$X_i$.}
\end{remark}

\begin{definition}
Let $A,F$ be formulas and $X \in {\cal P}$. The
formula $A[F/X]$ represents the result of substitution of $F$ to each
occurrence of $X$.
\end{definition}

We have the following result.

\begin{theorem}
Let $\G \cup \{A,F\}$ be a set of formulas, $X_m$ a minimal variable,
$X_i$ an intuitionistic variable, $X_c$ a classical variable, $F_i$ an
intuitionistic formula, and $F_c$ a classical formula. If $\G \v_{pml}
A$, then $\G[F/X_m] \v_{pml} A[F/X_m]$, $\G[F_i/X_i] \v_{pml}
A[F_i/X_i]$ and $\G[F_c/X_c] \v_{pml} A[F_c/X_c]$.
\end{theorem}

\begin{proof*}
By induction on the proof of $\G \v_{pml} A$.
\end{proof*}

\section{A semantics for PML}

Now we are ready for a definition of Kripke semantics for PML.

\begin{definition}
 A \emph{mixed Kripke model} is a triple ${\cal
K}=(K,{\pp},{\realise})$, where $(K,{\pp})$ is an inhabited, partially
ordered set (poset), and ${\realise}$ a binary relation on $K \times
{\cal P}$ such that:
\begin{janclenum}{\arabic}{)}
\item For all $\X \in {\cal P}$, if $\a \realise \X$ and $\b \pg
\a$, then $\b \realise \X$.

\item If $\a \realise \perp$, then, for all classical or
intuitionistic variable $X_s$, $\a \realise X_s$.

\item If $\a \realise X_c$ and, $\a \not \realise \perp$, then for
each $\b \in K$: $\b \realise X_c$.
\end{janclenum}
The relation $\realise$ is then extended to logically compound formulas by the following clauses:
\begin{itemize}

\item $\a \realise A \et B $ iff $\a \realise A$ and $\a \realise
B$.

\item  $\a \realise A \ou B $ iff $\a \realise A$ or $\a \realise
B $.

\item $\a \realise A \f B $ iff for all $\b \pg \a$ , if $\b
\realise A$, then $\b \realise B$.
\end{itemize}
\end{definition}

\begin{lemma}
For all formulas we have monotonicity: for all $\a,\b \in K$ ($\a
\realise A$ and $\b \geq \a$ implies $\b \realise A$).
\end{lemma}

\begin{proof*}
By formula induction.
\end{proof*}

\begin{definition}
A formula $A$ is \emph{valid} in a mixed Kripke model ${\cal
K}=(K,{\pp},{\realise})$ iff for all $\a \in K$, $\a \realise A$;
notation ${\cal K} \realise A$. If $\G$ is a set of formulas, we
say that $\G \realise A$ iff in each  mixed model ${\cal K}$ such
that: if for all $B \in \G$, ${\cal K} \realise B$, then also
${\cal K} \realise A$.
\end{definition}

\begin{remark} {To check if ${\cal K} \realise A$ it is enough to limit ${\cal K}$ to
the variables of $A$.}
\end{remark}

We have the following lemmas.
\begin{lemma} \label{kripkeint}
Let $A$ be an intuitionistic formula and ${\cal K}$ a  mixed Kripke
model. We have ${\cal K} \realise \perp \f A$.
\end{lemma}

\begin{proof*}
By induction on the complexity of $A$.
\end{proof*}

\begin{lemma}  \label{kripkecla}
Let $A$ be a classical formula and ${\cal K}$ a  mixed Kripke model.
We have ${\cal K} \realise \non \non A \f A$.
\end{lemma}

\begin{proof*}
We first prove, by induction, that if $B$ is a classical formula,
$\b \in K$ and $\b \realise B$, then, for each $\g \in K$, $\g
\realise B$. Let $\a \in K$ such that $\a \realise \non \non A$.
We may assume $\a \not \realise \perp$. Therefore $\a \not
\realise \non A$ and thus there is $\b \geq \a$ such that $\b
\realise A$. We deduce $\a \realise A$.
\end{proof*}

We can deduce the soundness theorem for PML.
\begin{theorem}
Let $\G \cup \{ A \}$ be a set of formulas.
If $\G \v_{pml} A$, then $\G \realise A$.
\end{theorem}

\begin{proof*}
The proof is by induction on derivation of $\G \v_{pml} A$ and we
use Lemmas \ref{kripkeint} and \ref{kripkecla}.
\end{proof*}

We present now a completeness proof for PML.
\begin{definition}
A set of formulas $\D$ is said to be \emph{saturated} iff: if $\D
\v_{pml} C \ou D$, then $C \in \D$ or $D \in \D$.
\end{definition}

\begin{remark} {A saturated set of formulas $\D$ is closed by deduction. Indeed, if $\D
\v_{pml} B$, then $\D \v_{pml} B \ou B$, thus $B \in \D$.}
\end{remark}

\begin{lemma}
If $\G \not \v_{pml}  A$, then there is a saturated set $\G_{\omega}$ such that $\G \inclus \G_{\omega}$ and
$\G_{\omega} \not \v_{pml}  A$.
\end{lemma}

\begin{proof*}
Same proof as the corresponding lemma in intuitionistic logic
\cite{DNR:2001,DAL:1994}.
\end{proof*}

\begin{definition}
Let $\G_0$ be any saturated set of formulas. Then we define ${\cal K}
= (K , \inclus , \realise)$ such that $K = \{ \D$ / $\D$ saturated sets
and $\G_0 \inclus \D \}$, and, for each $\X \in {\cal P}$: $\D
\realise \X$ iff $\X \in \D$.
\end{definition}

\begin{lemma} \label{kmodel}
${\cal K}$ is a  mixed Kripke model.
\end{lemma}

\begin{proof*} We must prove the three needed conditions:
\begin{enumerate}
\item Trivial.

\item If $\D \realise \perp$, then $\D \v_{pml} \perp$, thus $\D
\v_{pml} X_i$ and $\D \v_{pml} X_c$, i.e.\   $\D \realise X_i$ and
$\D \realise X_c$.

\item Let $\D \realise X_c$, $\D \not \realise \perp$, and $\D'
\not \realise \perp$. We have $\G_0 \v_{pml} X_c \ou \neg X_c$,
then $\G_0 \v_{pml} X_c$ or $\G_0 \v_{pml} \neg X_c$. Since $\G_0
\inclus \D$ and $\G_0 \inclus \D'$, we have $\G_0 \realise X_c$
and $\D' \realise X_c$.
\end{enumerate}
\end{proof*}

\begin{lemma}
For all $\D \in K$ and each formula $B$, $\D \realise B$ iff $B \in \D$.
\end{lemma}

\begin{proof*}
By induction on the complexity of $B$.
\end{proof*}

\begin{theorem}\label{comp}
Let $\G \cup \{ A \}$ be a set of formulas.
If $\G \realise A$, then $\G \v_{pml} A$.
\end{theorem}

\begin{proof*}
Suppose $\G \not \v_{pml}  A$, and let $\G_0$ be a saturated
extension of $\G$ such that $A \not \in \G_0$.  By  the last
construction there is a  mixed Kripke model ${\cal K} = (K ,
\inclus , \realise)$ and $\a \in K$ such that for all $B$: $\a
\realise B$ iff $B \in \G_0$. In particular, $\a \realise B$ for
$B \in \G$ and $\a \not \realise A$. Hence $\G \not \realise A$.
\end{proof*}

We also have the following results.
\begin{theorem}
\begin{enumerate}
\item The system {\rm PML} has the finite  mixed Kripke model
property.

\item The system {\rm PML} is decidable.
\end{enumerate}
\end{theorem}

\begin{proof*}
Same proof as the corresponding result in intuitionistic logic
\cite{DNR:2001,DAL:1994}.
\end{proof*}

\section{Properties of PML}

In this section we prove the principal result of the paper
(Theorems \ref{imp} and \ref{imp'}): ``To be derivable in the
system using only classical (resp.\ intuitionistic, minimal) rules
a mixed formula must contain at least a classical (resp.\
intuitionistic, minimal) variable''. This result is easily shown
if the system PML has some sub-formula property. However usually
such a property is a direct consequence of the cut-elimination
theorem which is difficult to show here because we cannot code the
disjunctive formulas (indeed the formula $\neg (\neg A \et \neg B)
\f A \ou B$ is not derivable) and eliminate the classical cuts.

\begin{definition}
\begin{janclenum}{(\arabic}{)}
\item An \emph{intuitionistic mixed Kripke model} (resp.\ a
\emph{minimal  mixed Kripke model}) is a mixed Kripke model
restricted on the formulas built on the set ${\cal P}_{(i)} =
{\cal V}_m \cup {\cal V}_i \cup \{ \perp \}$ (resp.\ the formulas
built on the set ${\cal P}_{(m)} = {\cal V}_m \cup \{ \perp \}$).

\item We write $\G \v_{(i)} A$ if $\G \v A$ is derivable without
using the rule $(\perp_c)$ and $\G \v_{(m)} A$ if $\G \v A$ is
derivable without using the rules $(\perp_i)$ and $(\perp_c)$.
\end{janclenum}
\end{definition}

We have the following results:
\begin{theorem} \label{comp'}
\begin{enumerate}
\item Let $\G \cup \{ A \}$ be a set of formulas without classical
variables. $\G \v_{(i)} A$ iff for all intuitionistic  mixed
Kripke model ${\cal K}$: ${\cal K} \realise \G$ implies ${\cal K}
\realise A$.

\item Let $\G \cup \{ A \}$ be a set of formulas without classical
and intuitionistic variables.  $\G \v_{(m)} A$ iff for all minimal
mixed Kripke model ${\cal K}$: ${\cal K} \realise \G$ implies
${\cal K} \realise A$.
\end{enumerate}
\end{theorem}

\begin{proof*}
In the proof of Theorem \ref{comp}, we use
the derivation rules to prove Lemma \ref{kmodel}.
\end{proof*}

\begin{definition}
For each  mixed Kripke model ${\cal K}$ we define the
intuitionistic (resp.\ the minimal)  mixed Kripke model ${\cal
K}_{(i)}$ (resp.\ ${\cal K}_{(m)}$) as being ${\cal K}$ restricted
on the set ${\cal P}_{(i)}$ (resp.\ ${\cal P}_{(m)}$). By
definition, it is clear that each intuitionistic  mixed Kripke
model (resp.\ minimal mixed Kripke model) can be seen as a ${\cal
K}_{(i)}$ (resp.\ a ${\cal K}_{(m)}$) for a  mixed Kripke model
${\cal K}$.
\end{definition}

\begin{lemma} \label{rest}
\begin{enumerate}
\item Let $A$ be a formula without classical variables. We have
${\cal K} \realise A$ iff ${\cal K}_{(i)} \realise A$.

\item Let $A$ be a formula without classical and intuitionistic
variables.  We have ${\cal K} \realise A$ iff ${\cal K}_{(m)}
\realise A$.
\end{enumerate}
\end{lemma}

\begin{proof*}
By induction on the complexity of $A$.
\end{proof*}

The following theorem is now an easy corollary.
\begin{theorem} \label{imp}
\begin{enumerate}
\item Let $\G \cup \{ A \}$ be a set of formulas without classical
variables. We have $\G \v_{pml} A$ iff $\G \v_{(i)} A$.

\item Let $\G \cup \{ A \}$ be a set of formulas without classical
and intuitionistic variables.  We have $\G \v_{pml} A$ iff $\G
\v_{(m)} A$.
\end{enumerate}
\end{theorem}

\begin{proof*}
\begin{enumerate}
\item If $\G \v_{pml} A$, then for all  mixed Kripke model ${\cal
K}$: ${\cal K} \realise \G$ implies ${\cal K} \realise A$, thus,
by Lemma \ref{rest}, for all intuitionistic mixed Kripke model
${\cal K}_{(i)}$: ${\cal K}_{(i)} \realise \G$ implies ${\cal
K}_{(i)} \realise A$. Therefore, by Theorem \ref{comp'}, $\G
\v_{(i)} A$.\label{enum:ref1}

\item Same proof as \ref{enum:ref1}).
\end{enumerate}
\vspace*{-1em}
\end{proof*}

\begin{definition}
We write $\G \v_{(i')} A$ if $\G \v A$ is derivable without using the rule
$(\perp_i)$.
\end{definition}

\begin{theorem} \label{imp'}
 Let $\G \cup \{ A \}$ be a set of formulas without intuitionistic variables.
$\G \v_{pml} A$ iff $\G \v_{(i')} A$.
\end{theorem}

\begin{proof*}
Same proof as Theorem \ref{imp}.
\end{proof*}

The proof of Theorem \ref{imp} is not constructive. We will try to make
a syntactical and constructive proof of this result (Corollary
\ref{imp-part}) but for a subsystem of PML.

\begin{definition}
Let ${\cal V'}_m$ be a countable subset of ${\cal V}_m$, and ${\tt m}$
be a bijective mapping between ${\cal V}_i$ and ${\cal V}'_m$. For all
formulas which do not contain classical variables the translation
$^{\tt m}$ is defined inductively by: $\perp^{\tt m} = \perp$,
${X_m}^{\tt m} = X_m$, ${X_i}^{\tt m} = \neg \neg {\tt m}(X_i)$ and
$(A \diamond B)^{\tt m} = A^{\tt m} \diamond B^{\tt m}$ if $\diamond
\in \{ \et , \ou , \f \}$.
\end{definition}

\begin{lemma}
Let $A$ be an intuitionistic formula.  $\v_{(m)} \perp \f A^{\tt m}$.
\end{lemma}

\begin{proof*}
By induction on $A$.
\end{proof*}

\begin{theorem}\label{tradm}
Let $\G \cup \{ A \}$ be a set of formulas without classical variables.
If $\G \v_{(i)} A$, then $\G^{\tt m} \v_{(m)} A^{\tt m}$.
\end{theorem}

\begin{proof*}
By induction on $\G \v_{(i)} A$.
\end{proof*}

\begin{corollary} \label{mincons}
Let $\G \cup \{ A \}$ be a set of formulas without classical and
intuitionistic variables. We have $\G \v_{(i)} A$ iff $\G \v_{(m)}
A$.
\end{corollary}

\begin{proof*}
By Theorem \ref{tradm}.
\end{proof*}

This method cannot be extended to get a syntactical proof of
Theorem \ref{imp}. We restrict our study to a subsystem of PML.

\begin{definition}
We denote by {\rm PML}$^{\ou}$ the system {\rm PML} with this
restriction on the rule $(\ou_E)$: if $A_1 \ou A_2$ is a classical
formula, then $B$ is also a classical formula.  We denote $\G
\v^{\ou}A$, if $A$ is derivable by $\G$ in {\rm PML}$^{\ou}$.
\end{definition}

\begin{remark}
The following derivation cannot be done in the system PML$^{\ou}$.
\begin{center}
{\tiny
$\F{\F{\F{
\F{\vdots}{\v X_c \ou \non X_c}
\quad
\F{
\F{}{X_c \v X_c} \quad
\F{}{X_c \f X_m \v X_c \f X_m}}
{X_c , X_c \f X_m \v  X_m}
\quad
\F{
\F{}{\neg X_c \v \neg X_c} \quad
\F{}{\neg X_c \f X_m \v \neg X_c \f X_m}}
{\neg X_c , \neg X_c \f X_m \v  X_m}}
{X_c \f X_m , \neg X_c \f X_m \v  X_m}}
{X_c \f X_m \v (\neg X_c \f X_m) \f  X_m}}
{\v (X_c \f X_m) \f ((\neg X_c \f X_m) \f  X_m)}$}
\end{center}
\end{remark}

\begin{definition}
Let ${\cal V}'_i$ be a countable subset of ${\cal V}_i$, and ${\tt i}$
be a bijective mapping between ${\cal V}_c$ and ${\cal V}'_i$. For all
formulas of {\rm PML} the translation $^{\tt i}$ is defined
inductively by: $\perp^{\tt i} = \perp$, ${X_m}^{\tt i} = X_m$,
${X_i}^{\tt i} = X_i$, ${X_c}^{\tt i} = \neg \neg {\tt i}(X_c)$, $(A
\diamond B)^{\tt i} = A^{\tt i} \diamond B^{\tt i}$ if $\diamond \in
\{ \et , \f \}$, and $(A \ou B)^{\tt i} = \neg \neg (A^{\tt i} \ou
B^{\tt i})$.
\end{definition}

\begin{lemma} \label{tradi}
Let $A$ be a classical formula.  We have $\v_{(i)} \neg \neg
A^{\tt i} \f A^{\tt i}$.
\end{lemma}

\begin{proof*}
By induction on $A$.
\end{proof*}

\begin{theorem}\label{intcons}
Let $\G \cup \{ A \}$ be a set of formulas.
If $\G \v^{\ou} A$, then $\G^{\tt i} \v_{(i)} A^{\tt i}$.
\end{theorem}

\begin{proof*}
By induction on $\G \v^{\ou} A$. We use Lemma \ref{tradi}
for the rules $(\perp_c)$ and $(\ou_E)$.
\end{proof*}

We can then deduce:
\begin{corollary} \label{imp-part}
\begin{enumerate}
\item Let $\G \cup \{ A \}$ be a set of formulas without classical
variables. If $\G \v^{\ou} A$, then $\G \v_{(i)}
A$.\label{enum:sqf1}

\item Let $\G \cup \{ A \}$ be a set of formulas without classical
and intuitionistic variables.  If $\G \v^{\ou} A$, then $\G
\v_{(m)} A$.\label{enum:sqf2}
\end{enumerate}
\end{corollary}

\begin{proof*}
\ref{enum:sqf1}) by Theorem \ref{intcons}, and \ref{enum:sqf2}) by
Corollary \ref{mincons}.
\end{proof*}

\section{Labels}

We establish in this section relations between PML and minimal,
intuitionistic and classical logics. If $A$ is a derivable formula
of ordinary propositional classical logic, we can label the
propositional variables of $A$ by $m$, $i$ or $c$ in order to
obtain a derivable formula in PML. It is clear that such a
labelling is not unique. We give in this section algorithms in
order to give ``minimal'' labels of classical propositional
formulas (Theorem \ref{decform}) and classical propositional
derivations (Theorem \ref{decdem}). We also define the notion of
``good'' derivation for a propositional classical formula
(Definition \ref{good}).

\begin{definition}
\begin{janclenum}{(\arabic}{)}
\item Let ${\cal V} = \{ X , Y , Z , ... \}$ be a countable set of
propositional variables. We suppose that ${\cal V}_m$ (resp.\
${\cal V}_i$, ${\cal V}_c$) are obtained by indexing the variables
of ${\cal V}$. Using ${\cal V} \cup \{\perp\}$ we define, as
usually, the minimal, intuitionistic, and classical logic denoted
respectively by {\rm PLM}, {\rm PLI} and {\rm PLC}. We use as
abbreviations $\v_m$, $\v_i$, $\v_c$ for derivability in {\rm
PLM}, {\rm PLI}, {\rm PLC} respectively. A formula built on ${\cal
V} \cup \{ \perp \}$ is called \emph{ordinary formula}.

\item A \emph{label} is a function $l : {\cal V} \f {\cal P}$ such
that $l(X) \in \{X_m , X_i , X_c \}$. A label $l$ is extended to
logical formulas by the following clauses: $l(\perp) = \perp$ and
$l(A \diamond B) = l(A) \diamond l(B)$ if $\diamond \in \{ \et ,
\ou , \f \}$.

\item We define on ${\cal V}_m \cup {\cal V}_i \cup {\cal V}_c$ a
binary relation $<$ as follows: for all $X \in {\cal V}$, $X_m <
X_i <X_c$. We define on labels a binary relation $<$ as follows:
$l < l'$ iff (1) for all variable $X \in {\cal V}$, $l(X) \leq
l'(X)$ and (2) there is a $X \in {\cal V}$ such that $l(X) <
l'(X)$.

\item Let $l_m$ (resp.\ $l_i$, $l_c$) be the label defined by: for
all $X \in {\cal V}$, $l_m(X) = X_m$ (resp.\ $l_i(X) = X_i$,
$l_c(X) = X_c$).
\end{janclenum}
\end{definition}

The following result means that PML contains disjoint copies of
systems PLM, PLI and PLC.
\begin{theorem}\label{dec1}
Let $\G \cup \{ A \}$ be a set of ordinary
formulas. We have: $\G \v_m A$ iff $l_m(\G) \v_{(m)} l_m(A)$,
$\G \v_i A$ iff $l_i(\G) \v_{(i)} l_i(A)$
and $\G \v_c A$ iff $l_c(\G) \v_{pml} l_c(A)$.
\end{theorem}

\begin{proof*}
Easy.
\end{proof*}

\begin{definition}
Let $A$ be an ordinary formula such that $\v_c A$. A \emph{label
for $A$} is a label $l$ such that $\v_{pml} l(A)$ and for every
variable $X$ which does not appear in $A$, $l(X) = X_m$.
\end{definition}

\begin{remark}
Let $A$ be an ordinary formula such that $\v_c A$. By Theorem
\ref{dec1}, $l_c$ is a label for $A$.
\end{remark}

\begin{definition}
Let $A$ be an ordinary formula such that $\v_c A$.  A
\emph{minimal label for $A$} is a label $l$ for $A$ such that: if
$l' \leq l$ is a label for $A$, then $l' = l$.
\end{definition}

\begin{theorem} \label{decform}
Let $A$ be an ordinary formula such that $\v_c A$. $A$ has a
minimal label.
\end{theorem}

\begin{proof*}
Since PML is decidable we try all possible labels for
$A$.
\end{proof*}

\begin{example}\label{dece}
Let $b$ the label defined by: $b(X) = X_c$, $b(Y) = Y_i$, and
for every $Z \not = X$ and $Y$, $b(Z) = Z_m$. It is easy to check that
$b$ is the unique minimal label for the ordinary formula $(Z \f
X) \ou (X \f Y)$. The minimal label for an ordinary formula is not
unique.  Let $A = (X \f Y) \ou (Y \f X)$ and $l,l'$ such
that $l(X)=X_c$, $l(Y)=Y_i$, $l'(X)=X_i$ and $l'(Y)=Y_c$. It is easy
to check that $l$ and $l'$ are two minimal labels for $A$ but
they are not comparable.
\end{example}

\begin{definition}
Let ${\cal D}$ be a derivation in {\rm PLC}. A \emph{label for
${\cal D}$} is a label $l$ such that: (1) for every variable $X$
which does not appear in ${\cal D}$, $l(X) = X_m$ and (2) by
extending $l$ on ${\cal D}$ we obtain a derivation in {\rm PML}. A
\emph{minimal label for ${\cal D}$} is a label $l$ for ${\cal D}$
such that: if $l' \leq l$ is a label for ${\cal D}$, then $l' =
l$.
\end{definition}

\begin{remark}  {$l_m$ (resp.\ $l_i$, $l_c$) is a label for all
 derivation in {\rm PLM} (resp.\ {\rm PLI}, {\rm PLC}).}
\end{remark}

\begin{definition}
Let $l_1,...,l_n$ be labels. We define a new label
$sup(l_1,...,l_n)$ as follows: for every $X \in {\cal V}$,
$sup(l_1,...,l_n)(X) = sup(l_1(X),...,l_n(X))$.
\end{definition}

\begin{theorem} \label{decdem}
Let ${\cal D}$ be a derivation in {\rm PLC}. The derivation ${\cal D}$ has a
unique minimal label.
\end{theorem}

\begin{proof*}
We define the unique minimal label $l_{\cal D}$ by induction on
${\cal D}$.
\begin{enumerate}
\item If ${\cal D}$ is $(Ax)$, then $l_{\cal D} = l_m$.

\item If the last rule used in ${\cal D}$ is
    \begin{itemize}
    \item $(W)$, $(\et_E)$, $(\ou_I)$, or $(\f_I)$, then $l_{\cal D} =
    l_{{\cal D}_1}$.

    \item  $(\et_I)$, or $(\f_E)$, then $l_{\cal D} = sup(l_{{\cal
    D}_1},l_{{\cal D}_2})$.

    \item $(\ou_E)$, then $l_{\cal D} = sup(l_{{\cal
    D}_1},l_{{\cal D}_2},l_{{\cal D}_3})$.

    \item $(\perp_i)$, then $l_{\cal D} = l \circ l_{{\cal
    D}_1}$, where
    \item[] $l(l_{{\cal D}_1}(X)) =
    \begin{cases}
    X_i & \text{if }X \in var(A) \text{ and }l_{{\cal D}_1}(X) \not =
    X_c\\
    X_c & \text{if }X \in var(A) \text{ and }l_{{\cal D}_1}(X) = X_c\\
    l_{{\cal D}_1}(X) & \text{otherwise }
    \end{cases}$

    \item $(\perp_c)$, then $l_{\cal D} = l \circ l_{{\cal D}_1}$,
    where $l(l_{{\cal D}_1}(X)) =
    \begin{cases}
    X_c     &   \text{if }X \in var(A)\\
    l_{{\cal D}_1}(X) & \text{otherwise }
    \end{cases}$
    \end{itemize}
\end{enumerate}
\end{proof*}

\begin{example} \label{ex}
It is easy to check that the label $b$ of the Example \ref{dece}
is the minimal label for the following derivation:
\begin{center}
{\footnotesize $\F{\F{\vdots}{\v X \ou \non X} \quad
\F{\F{\F{\F{}{X \v X}}{X , Z \v X}}{X \v Z \f X}}{X
 \v (Z \f X) \ou (X \f Y)}
 \quad
\F{\F{\F{\F{\F{}{X \v X} \quad \F{}{\neg X \v \neg X}}{X , \neg X \v \perp}}{X, \neg X \v Y}}{
\neg X \v  X \f Y}}{\neg X \v (Z \f X) \ou (X \f Y)}}
{\v (Z \f X) \ou (X \f Y)}$}
\end{center}
\end{example}

\begin{definition} \label{good}
Let $A$ be an ordinary formula such that $\v_c A$. A \emph{good
derivation} for $A$ is a derivation ${\cal D}$ of $A$ in {\rm PLC}
such that $l_{\cal D}$ is a minimal label for $A$. Intuitively, a
good derivation of a formula $A$ is a derivation of $A$ with
minimal use of the absurdity rules.
\end{definition}

\begin{theorem}
Let $A$ be an ordinary formula such that $\v_c A$. The formula $A$ has
a good derivation.
\end{theorem}

\begin{proof*} Let $l_A$ be a minimal label of $A$. Since we can
enumerate all derivable formulas, then we can find a derivation
${\cal D}$ ending with $l_A(A)$. The derivation obtained by
erasing the indexes in the derivation ${\cal D}$ is a good
derivation for $A$.
\end{proof*}

\begin{example}
The derivation of the Example \ref{ex} is a good derivation for
the formula $(Z \f X) \ou (X \f Y)$.
\end{example}

\section{Sequent calculus}

We describe below a sequent calculus version of PML. This sequent
calculus is non satisfactory because it does not satisfy  the
cut-elimination property (Theorem \ref{cut}).

\begin{definition}
In this section a sequent is of the form $\G \v' A; \D$ where $\G$
(resp.\ $\D$) is a finite set of formulas (resp.\ of classical
formulas) and $A$ is a formula. The rules of sequent calculus are
the following:

\begin{minipage}[t]{160pt}
$(Ax) \; \F{}{A \v' A;}$\\
\end{minipage}
$(Cut) \; \F{\G_1 ,A \v' B ; \D_1   \;\;\;  \G_2 \v' A ; \D_2 }
{ \G_1 , \G_2  \v' B ; \D_1 , \D_2}$\\
\begin{minipage}[t]{160pt}

\end{minipage}

\begin{minipage}[t]{160pt}
$(S_r) \; \F{\G \v'A ; \perp , \D } { \G \v'A ; \D}$
\end{minipage}
\begin{minipage}[t]{140pt}
$(S_l) \; \F{\G \v'A ; A , \D } { \G \v' A ; \D}$\\
\end{minipage}

\begin{minipage}[t]{200pt}
$(W_r) \; \F{\G \v' \perp ; \D \;\;\; A {\rm \; is \; an \; intuitionistic \; formula}    } { \G \v' A ; \D}$\\
\end{minipage}

\begin{minipage}[t]{160pt}
$(W_l) \; \F{\G \v' A ; \D } { \G , B \v' A ; \D}$\\
\end{minipage}

\begin{minipage}[t]{180pt}
$(W'_r) \;\F{\G \v' A ; \D  \;\;\; B {\rm \; is \; a \; classical \; formula} } { \G \v' A ; B, \D}$\\
\end{minipage}

\begin{minipage}[t]{190pt}
$(E) \; \F{\G \v'A ; B , \D \;\;\; A {\rm \; is \; a \; classical \; formula}  } { \G \v'B ; A , \D}$\\
\end{minipage}

\begin{minipage}[t]{155pt}
$(\et_r) \; \F{\G_1 \v' A_1 ; \D_1  \;\;\;  \G_2 \v' A_2 ; \D_2}
{ \G_1 , \G_2 \v' A_1 \et A_2 ; \D_1 , \D_2 }$\\
\end{minipage}
\begin{minipage}[t]{160pt}
$(\et_l) \; \F{\G , A_i \v' B ; \D  } { \G , A_1 \et A_2 \v' B ; \D }$ \\
\end{minipage}

\begin{minipage}[t]{155pt}
$(\ou_r) \;\F{\G \v' A_i ; \D} { \G \v' A_1 \ou A_2 ; \D }$
\end{minipage}
\begin{minipage}[t]{165pt}
$(\ou_l) \; \F{\G_1 , A_1 \v' B ; \D_1\;\;\;  \G_2 , A_2 \v' B ; \D_2}
{ \G_1 , \G_2, A_1 \ou A_2 \v' B ; \D_1 , \D_2}$\\
\end{minipage}

\begin{minipage}[t]{155pt}
$(\f_r) \; \F{\G , A_1 \v' A_2 ; \D} { \G \v' A_1 \f A_2 ; \D }$\\
\end{minipage}
\begin{minipage}[t]{160pt}
$(\f_l) \; \F{\G_1 \v' A_1 ; \D_1 \;\;\; \G_2 , A_2 \v' B ; \D_2}
{ \G_1 , \G_2 , A_1 \f A_2 \v' B ; \D_1 , \D_2 }$\\
\end{minipage}

We write $\G \v^{pml} A ; \D$ if there is a derivation ${\cal D}$
ending with the sequent $\G \v' A ; \D$.
\end{definition}

We wish to show $\G \v^{pml} A;$ iff $\G \v_{pml} A$.
\begin{lemma} \label{exSC}
\begin{enumerate}
\item If $A$ is an intuitionistic formula, then $\v^{pml}\perp \f
A;$.\label{enum:qss1}

\item If $B$ is a classical formula, then $\v^{pml} \non \non B \f
B;$.\label{enum:qss2}
\end{enumerate}
\end{lemma}

\begin{proof*} \ref{enum:qss1}) is easy. For \ref{enum:qss2}):
\begin{center}
{\footnotesize
 $\F{\F{\F{\F{\F{\F{\F{\F{}{B \v' B;}}{B \v' B;
\perp}}{B \v' \perp ; B}} {\v' \non B ; B} \quad \F{}{\perp \v'
\perp;}} {\non \non B \v' \perp; B}} {\non \non B \v' B; \perp}}
{\non \non B \v' B;}} {\v' \non \non B \f B;}$}
\end{center}
\vspace*{-1em}
\end{proof*}

\begin{theorem} \label{dirc}
Let $\G \cup \{ A \}$ be a set of formulas. If $\G \v_{pml} A$, then $\G \v^{pml} A;$.
\end{theorem}

\begin{proof*}
By induction on the proof of $\G \v_{pml} A$. We use the
cut rule and Lemma \ref{exSC}.
\end{proof*}

\begin{lemma} \label{exSC'}
If $A,B$ are classical formulas, then $\v_{pml} [(\non A \f A) \f A]
\et [(\non B \f A) \f (\non A \f B)]$.
\end{lemma}

\begin{proof*}
Easy.
\end{proof*}

\begin{definition}
Let $\neg \D$ indicate the negation of the formulas in $\D$.
\end{definition}

\begin{theorem} \label{recip}
Let $\G$ be a set of formulas, $\D$ a set of classical formulas, and
$A$ a formula. If $\G \v^{pml} A; \D$, then $\G , \neg \D \v_{pml} A$.
\end{theorem}

\begin{proof*}
By induction on the proof of $\G \v^{pml} A; \D$. We use
Lemma \ref{exSC'} for the rules $(E)$ and $(S_l)$.
\end{proof*}

We can then deduce:
\begin{corollary}
Let $\G \cup \{ A \}$ be a set of formulas.
We have $\G \v_{pml} A$ iff $\G \v^{pml} A;$.
\end{corollary}

\begin{proof*}
We use Theorems \ref{dirc} and \ref{recip}.
\end{proof*}

\begin{remark}
The usual process to eliminate cuts in the sequent calculus is not
valid for our system. For example, the elimination of cuts in the
following derivation needs the use of several non classical
formulas on the right.
\begin{center}
{\footnotesize $\F{\F{\F{\F{\F{\F{\F{\F{\F{}{X_c \v' X_c;}}{X_c
\v' X_c;\perp}}{X_c \v'\perp; X_c}}{\v'\neg X_c; X_c}}{\v'X_c \ou
\neg X_c; X_c}} {\v'X_c; X_c \ou \neg X_c}}{\v'X_c \ou \neg X_c;
X_c \ou \neg X_c}}{\v'X_c \ou \neg X_c;} \quad \F{\F{\F{}{X_c \v'
X_c;} \quad \F{}{X_m \v' X_m;}}{X_c , X_c \f X_m \v' X_m;} \quad
\F{\F{}{\neg X_c \v' \neg X_c;} \quad \F{}{X_m \v' X_m;}}{\neg X_c
, \neg X_c \f X_m \v' X_m;}} {X_c \ou \neg X_c , X_c \f X_m, \neg
X_c \f X_m \v' X_m; }} {X_c \f X_m, \neg X_c \f X_m \v' X_m;}$}
\end{center}
\end{remark}

\begin{theorem}\label{cut}
The {\rm PML} sequent calculus does not satisfy the
cut-elimination (even weak) property.
\end{theorem}

\begin{proof*} We prove that there is no normal derivation (i.e.\  without
cuts) for the sequent $ X_c \f X_m , \neg X_c \f X_m \v ^{pml} X_m ;$.
By using the following mixed Kripke model $\cal K$ = $(K , \pp ,
\realise)$ where $ K = \{ \alpha , \beta \}$, $\alpha \pp \beta$,
$\beta \realise X_m$, and $\beta \realise \faux$, we prove easily that
$X_c \f X_m \not \v ^{pml} X_m ;$, $X_c \f X_m \not \v ^{pml} \neg X_c
;$, $\neg X_c \f X_m \not \v ^{pml} X_c ;$, $\neg X_c \f X_m \not \v
^{pml} X_m ;$, $\not \v ^{pml} X_c ;$, and $\not \v ^{pml} \neg X_c
;$.  Let us take a minimal derivation of $ X_c \f X_m , \neg X_c \f
X_m \v ^{pml} X_m ;$ and look at the last used rule.
\begin{enumerate}
\item If it is the rule $(W_l)$, then $X_c \f X_m \v ^{pml} X_m ;$
or $\neg X_c \f X_m \v ^{pml} X_m ;$.

\item If it is the rule $(\f _l)$, then $\neg X_c \f X_m \v ^{pml}
X_c ;$ or $\v ^{pml} X_c ;$ or $X_c \f X_m \v ^{pml} \neg X_c ;$
or $ \v ^{pml} \neg X_c ;$.

\item If it is the rule $(S_r)$, then $ X_c \f X_m , \neg X_c \f
X_m \v ^{pml} X_m ; \faux $. We again look at the last rule used.
\begin{itemize}
    \item If it is the rule $(W_l)$, then $X_c \f X_m \v ^{pml} X_m ; \faux $
or $\neg X_c \f X_m \v ^{pml} X_m ; \faux $.

    \item If it is the rule $(\f _l)$, then $\neg X_c \f X_m \v ^{pml} X_c
;$ or $\neg X_c \f X_m \v ^{pml} X_c ; \faux$ or $\v ^{pml} X_c ;$
or $\v ^{pml} X_c ; \faux $ or $X_c \f X_m \v ^{pml} \neg X_c ;$
or $X_c \f X_m\v ^{pml} \neg X_c ; \faux$ or $ \v ^{pml} \neg X_c
;$ or $ \v ^{pml} \neg X_c ;\faux$.
\end{itemize}
\end{enumerate}
\end{proof*}

\begin{remark}
To get a normal derivation of the sequent $X_c \f X_m, \neg X_c \f
X_m \v' X_m;$,  we need more flexible rules. For example:
\begin{itemize}
\item allowing the use of the logical rules each formula on the
right;

\item allowing several occurrences of the same non classical
formula on the right.
\end{itemize}
Here is a derivation of sequent $X_c \f  X_m, \neg X_c \f  X_m \v'
 X_m$ without using the cut rule.
\begin{center}
{\footnotesize $\F{\F{\F{\F{\F{\F{}{X_c \v' X_c} \quad \F{}{X_m
\v' X_m}}{X_c , X_c \f  X_m \v' X_m}} {X_c , X_c \f  X_m \v' X_m ,
\perp}} {X_c \f  X_m \v' X_m , \neg X_c} \quad \F{}{X_m \v'  X_m}}
{X_c \f  X_m , \neg X_c \f X_m \v' X_m , X_m}} {X_c \f  X_m , \neg
X_c \f  X_m \v' X_m }$}
\end{center}
\end{remark}

\begin{oq*}
``Is it possible to eliminate cuts in
such a system?''
\end{oq*}

\acknowledgements{ We wish to thank Noël Bernard, René David,
François Pabion and Christophe Raffalli for helpful discussions.
We also thank Jamil Nour for his help in the writing of this
paper.}


\end{document}